\documentclass[12pt]{amsart}
\usepackage{amsfonts, latexsym, amssymb, amsgen, amsbsy, amstext, amsopn, amsmath}
\usepackage{color}

\newcommand{\E}{\mathbb{E}}
\newcommand{\G}{\mathbb{G}}

\newcommand{\N}{\mathbb{N}}

\newcommand{\R}{\mathbb{R}}

\newcommand{\cE}{\mathcal{E}}

\newcommand{\cL}{\mathcal{L}}

\newcommand{\cS}{\mathcal{S}}

\renewcommand{\span}{\mbox{\rm span}}

\newcommand{\ph}{\varphi}

\newcommand{\sm}{\setminus}

\newcommand{\res}{\mbox{\LARGE{$\llcorner$}}}
\newcommand{\lan}{\langle}
\newcommand{\ran}{\rangle}

\newcommand{\lra}{\longrightarrow}

\newcommand{\der}{\partial}

\newcommand{\Bx}{\mbox{\rm Box}}

\newtheorem{The}{Theorem}
\newtheorem{Lem}[The]{Lemma}

\newtheorem{Rem}[The]{Remark}
\newtheorem{Pro}[The]{Proposition}

\newtheorem{Cor}[The]{Corollary}

\begin{document}
\title{Measure of submanifolds in the Engel group}

\author{Enrico Le Donne}
\address{Yale University, USA}
\email{enrico.ledonne@yale.edu}

\author{Valentino Magnani}
\address{Valentino Magnani: Dipartimento di Matematica \\
Largo Pontecorvo 5 \\ I-56127, Pisa}
\email{magnani@dm.unipi.it}

\begin{abstract}
We find all intrinsic measures of $C^{1,1}$ smooth
submanifolds in the Engel group, showing that they are equivalent
to the corresponding $d$-dimensional spherical Hausdorff measure
restricted to the submanifold.
The integer $d$ is the degree of the submanifold. 
These results follow from a different approach to negligibility,
based on a blow-up technique.
\end{abstract}

\maketitle

\section{Introduction}

Computing the Hausdorff measure of submanifolds in stratified groups
with respect to the Carnot-Carath\'eodory distance is a rather natural question.
This may be considered as a first step to study several problems of
Geometric Measure Theory in stratified groups.

Nevertheless, this question has not yet an answer.
In 0.6~B of \cite{Gr1}, Gromov has given a general formula for
the Hausdorff dimension
of smooth submanifolds in equiregular Carnot-Carath\'eodory spaces
and in \cite{Mag8B} it is shown that this formula coincides with the degree of
the submanifold, recently introduced in \cite{MagVit}.
In the latter work, the authors find an integral formula for the spherical
Hausdorff measure of submanifolds in stratified groups under a suitable 
``negligibility condition''.
If $d$ is the degree of a submanifold, and $\cS^d$ is the 
spherical Hausdorff measure constructed with the Carnot-Carath\'eodory distance,
this condition requires that all points of the submanifold having pointwise degree 
less than $d$ must be $\cS^d$-negligible. This negligibility condition has been recently
obtained in all two step groups, \cite{Mag8B}, but it is still open in higher step groups.
We address the reader to \cite{MagVit}, \cite{Mag8A} and \cite{Mag8B}
for more information on this problem and its connections with the
present literature.

In this work, we prove the negligibility condition in the Engel group,
adopting a different approach with respect to the standard covering arguments.
Broadly speaking, we simply ``blow-up'' the points of the submanifold adopting 
the intrinsic dilations of the group and then apply a simple fact of Geometric 
Measure Theory, see Lemma~\ref{gmta}. Essentially, we prove that the assumptions
of this lemma hold in all the single cases that can occur.
In fact, joining all propositions
of Section~\ref{sectsurfaces} and \ref{sectcurves}, we have our main result.
\begin{The}\label{liminfty}
Let $\Sigma$ be a $p$-dimensional $C^{1,1}$ submanifold of degree $d$
in the Engel group, where $p=1,2$.
Then for every $x\in\Sigma$ with $d_\Sigma(x)<d$, there holds
\begin{equation}\label{rdinfty}
\lim_{r\to0}\frac{\mu_p(\Sigma\cap D_{x,r})}{r^d}=+\infty.
\end{equation}
\end{The}
We have denoted by $\mu_p$ the $p$-dimensional Riemannian surface measure 
induced on $\Sigma$ by a fixed left invariant Riemannian metric on the group.
The main feature of the previous theorem is that \eqref{rdinfty}
not only depends on the degree at $x$, but also on the ``behaviour'' of $\Sigma$
around $x$, that is expressed by the degree of $\Sigma$.
It is rather interesting to observe how the limit \eqref{rdinfty} in some cases 
requires the use of this ``global'' information.
This is the case of Proposition~\ref{deg3deg2}, where the fact that
$\Sigma$ has degree three implies a differential constraint, given by
the system of PDEs \eqref{deg3pdes}, that play a crucial role in the proof of \eqref{rdinfty}.

Theorem~\ref{liminfty} joined with Lemma~\ref{gmta} yields Theorems~\ref{neglsurf} and
\ref{neglcurve}, that correspond to the negligibility condition stated in (1.5) of
\cite{MagVit}. This condition gives (1.4) of \cite{MagVit}, namely, 
we have the following
\begin{Cor}\label{areatypeforengel}
Let $\Sigma$ be a $p$-dimensional $C^{1,1}$ submanifold of degree $d$ in the Engel group.
Then the following formula holds
\begin{equation}\label{integrfor}
\int_{\Sigma} \theta(\tau_{\Sigma}^d(x))\, d\cS^d(x)
=\int_{\Sigma}|\tau_{\Sigma}^d(x)|\,d\tilde\mu_p(x).
\end{equation}
\end{Cor}
The $p$-tangent vector $\tau_{\Sigma}^d(x)$ is the component of degree $d$
of the $p$-tangent vector $\tau_\Sigma(x)$ associated to the tangent space $T_x\Sigma$.
Its norm is computed with respect to the auxiliary Riemannian metric fixed on the group.
This metric also yields the surface measure $\tilde\mu_p$ induced on $\Sigma$,
see \cite{MagVit} for more details.
In the case $p=3$, the previous integral formula follows 
from Theorem~2.20 of \cite{Mag5}. In fact, Frobenius Theorem implies that
$C^{1,1}$ hypersurfaces in any stratified group must possess non-horizontal points,
hence they have degree equal to $Q$$-$1, where $Q$ is the Hausdorff dimension
of the group. Recall that non-horizontal points have been introduced in
in \cite{MagVit} and studied in \cite{Mag8A}. 
According to Proposition~3.2 of \cite{Mag8B}, 
the length of the horizontal normal $|n_H(x)|$ in Theorem~2.20 of \cite{Mag5}
is equal to the length 
$|\tau_{\Sigma}^d(x)|$. Of course, $4$-dimensional submanifolds of the
Engel group are just open subsets, for which it is trivial to observe that
their degree is exactly 7 and their 7-dimensional Hausdorff measure is clearly 
positive and finite on the intersection with bounded sets.
The metric factor $\theta(\tau_{\Sigma}^d(x))$ is uniformly
bounded from above and from below, then $\cS^d\res\Sigma$ is equivalent to
the intrinsic measure $|\tau_{\Sigma}^d(x)|\,\tilde\mu_p\res\Sigma$,
introduced in \cite{MagVit}.
In case it is possible to find a distance that yields a constant metric factor,
then up to a geometric constant, we obtain
\begin{equation}\label{integrfor1}
\cS^d_\G(\Sigma)=\int_{\Sigma}|\tau_{\Sigma}^d(x)|\,d\tilde\mu_p(x).
\end{equation}
As an immediate consequence, the degree of a $C^{1,1}$ submanifold in the
Engel group equals its Hausdorff dimension, since points of maximum degree
form an open subset of the submanifold. In Remark~\ref{gromovproblemEngel}
we point out how our results are also related to the
Gromov's dimension comparison problem recently raised in \cite{BTW}.


%
%
%
%
%
\section{Basic definitions and standard results}\label{basdef}
%
%
%
%
%

The Engel group $\E$ is a connected, simply connected stratified group,
whose Lie algebra satisfies the direct decomposition
\[
\mathcal E=V_1\oplus V_2\oplus V_3
\]
and there exists a basis $(X_1,X_2,X_3,X_4)$ of $\mathcal E$, 
such that the only nontrivial brackets are
\[
[X_1,X_2]=X_3\quad\mbox{and}\quad [X_1,X_3]=X_4,
\]
where $V_1=\span\{X_1,X_2\},$ $V_2=\span\{X_3\}$ and $V_3=\span\{X_4\}$.
We represent the Engel group $\E$ by $\R^4$ equipped with the vector fields
\begin{eqnarray}\label{leftEng}
X_1=\der_1,\quad X_2=\der_2+x_1\der_3+\frac{x_1^2}{2}\der_4,\quad X_3=\der_3+x_1\der_4,
\quad X_4=\der_4,
\end{eqnarray}
where the associated exponential mapping builds the group operation in
$\R^4$ that makes it isomorphic to the abstract Engel group.
The intrinsic dilations $\delta_r:\R^4\lra\R^4$
are given by $\delta_r(x)=(rx_1,rx_2,r^2x_3,r^3x_4)$, with $r>0$.
This is a one parameter family of group automorphisms, since
\[
(\delta_r)_*\big(X_1\big)=rX_1,\quad(\delta_r)_*\big(X_2\big)=rX_2,
\quad(\delta_r)_*\big(X_3\big)=r^2X_3,\quad (\delta_r)_*\big(X_4\big)=r^3X_4,
\]
as one can check from direct computation.
We fix a left invariant Riemannian metric in $\R^4$ that makes $X_j$'s orthonormal. 
The Carnot-Carath\'eodory distance associated to
$\span\{X_1,X_2\}$ along with the fixed left invaraint metric on $\R^4$
yields a homogeneous distance. More generally, we will consider an
arbitrary homogeneous distance $\rho$ on $\R^4$,
namely, a continuous, left invariant distance that satisfies
\[
\rho(\delta_rx,\delta_ry)=r\,\rho(x,y)\quad\mbox{for every}\quad x,y\in\R^4,
\quad r>0.
\]
In the sequel, the abstract Engel group $\E$ will be identified with
$\R^4$, equipped with left invariant vector fields \eqref{leftEng},
distance $\rho$ and dilations $\delta_r$. The explicit formula for the
group operation in $\R^4$ will be not needed.

Our arguments are based on the following elementary fact of
Geometric Measure Theory, see for instance 2.10.19 of \cite{Fed}.
\begin{Lem}
Let $X$ be a metric space, let $\mu$ be a Borel measure on $X$
and let $\{V_i\}_{i\in\N}$ be an open covering of $X$
such that $\mu(V_i)<\infty$. Let $Z\subset X$ be a Borel set and
suppose that
\[
\limsup_{r\to0^+}r^{-a}\mu(D_{x,r})\geq\kappa>0
\]
whenever $x\in Z$, where $a>0$. Then $\mu(Z)\geq\kappa\;\cS^a(Z)$.
\end{Lem}
We have denoted by $\cS^a$ the $a$-dimensional spherical Hausdorff measure
constructed with the size function $\zeta_a(D_{x,r})=r^a$ and
$D_{x,r}$ is the closed ball of center $x$ and radius $r$.
From the previous lemma, we get the straightforward
\begin{Lem}\label{gmta}
Let $\Sigma$ be $k$-dimensional $C^{1,1}$ submanifold of $\E$ and
let $\mu_k$ be the left invariant Riemannian measure of $\E$
restricted to $\Sigma$. If $Z$ is a Borel set of $\Sigma$ such that
$\limsup_{r\to0^+}r^{-a}\mu_k(D_{z,r})=+\infty$, whenever
$z\in Z$, then $\cS^a(Z)=0$.
\end{Lem}

%
%
%
%
%
%
%
%
%
\section{Degree of submanifolds in the Engel group}
%
%
%
%
%
%

The degree of a 2-vector
$\tau=\sum_{1\leq i<j\leq4} \tau_{ij}\; X_i\wedge X_j\in\Lambda_2(\cE)$
is given by
\[
\deg(\tau)=\max\{d_i+d_j\mid \tau_{ij}\neq0\}
\]
where $d_i$ is the degree of $X_i$, hence $d_1=d_2=1$, $d_3=2$ and $d_4=3$.
Analogously, the degree of a vector $\tau=\sum_{i=1}^4 \tau_i\; X_i\in\cE$
is given by $\deg(\tau)=\max\{d_i\mid \tau_i\neq0\}$.
Then we define the pointwise degree at $x$ of a 
$p$-dimensional submanifold $\Sigma$ in $\E$ as
\[
d_\Sigma(x)=\deg\big(\tau_\Sigma(x)\big),
\]
where $\tau_\Sigma(x)$ is the $p$-tangent vector of $\Sigma$ at $x\in\Sigma$,
$p=1,2$.
If $\Sigma$ is a submanifold of $\E$ we define its degree $d(\Sigma)$ as
the integer $\max_{x\in\Sigma}d_\Sigma(x)$, see \cite{MagVit} for more details
in the general case of stratified groups. 
Let $U$ be an open subset of $\R^2$ and let $\phi:U\lra\R^4$ be a
$C^1$ immersion. According to computations in Section~4 of \cite{MagVit}, we have
\begin{eqnarray}\label{formula}
&&\phi_{u_1}\wedge\phi_{u_2}=
\phi_u^{12} X_1\wedge X_2+\left(\phi_u^{13}-\phi_1\phi_u^{12}\right)
X_1\wedge X_3+\phi_u^{23}  X_2\wedge X_3\\
&&+\left(\phi_u^{14}-\phi_1\,\phi_u^{13}+\frac{(\phi_1)^2}{2}\,
\phi_u^{12}\right)X_1\wedge X_4
+\left(\phi_u^{24}-\phi_1\,\phi_u^{23}\right) X_2\wedge X_4 \nonumber\\
&&
+\left(\phi_u^{34}+\frac{(\phi_1)^2}{2}\phi_u^{23}
-\phi_1\phi_u^{24}\right) X_3\wedge X_4\,,\nonumber
\end{eqnarray}
where we have defined
\[
\phi=(\phi_1,\phi_2,\phi_3,\phi_4)\quad\mbox{and}\quad
\phi^{ij}_u=\det\left(\begin{array}{cc}\phi^i_{u_1} & \phi^i_{u_2} \\
\phi^j_{u_1} & \phi^j_{u_2}\end{array}\right).
\]
It is also understood that $X_i\wedge X_j$ in the previous formula
are evaluated at the point $\phi(u)$.
Thus, if $\phi$ locally parametrizes a submanifold $\Sigma$, according to the
notion of pointwise degree, we have that
\begin{equation}\label{ciju}
d_\Sigma(\phi(u))=
\left\{\begin{array}{ll} 5 & \mbox{if $c_{34}(u)\neq0$} \\ 
4 & \mbox{if $|c_{14}(u)|+|c_{24}(u)|>0$ and $c_{34}(u)=0$} \\ 
3 & \mbox{if $|c_{13}(u)|+|c_{23}(u)|>0$ and $c_{34}(u)=c_{14}(u)=c_{24}(u)=0$} \\
2 & \mbox{if $c_{34}(u)=c_{14}(u)=c_{24}(u)=c_{13}(u)=c_{23}(u)=0$} \end{array}\right.\,,
\end{equation}
where we have set 
\begin{eqnarray}\label{formulac}
\phi_{u_1}\wedge\phi_{u_2}=\sum_{1\leq i<j\leq4} c_{ij}(u)\; X_i\wedge X_j\,.
\end{eqnarray}
\begin{Rem}\label{gromovproblemEngel}{\rm
By definition of degree, in the Engel group one easily notices that 
all possible degrees of $C^{1,1}$ surfaces can only be 3,4 or 5. 
Degree two is not possible due to the Frobenius Theorem. 
Of course, curves can only have degrees 1,2 or 3
and again Frobenius Theorem implies that hypersurfaces can only have degree 6.
Thus, these are all possible Hausdorff dimensions of $C^{1,1}$ smooth
submanifolds in the Engel group and formula \eqref{integrfor} holds for them.
This answers Problem~1.1 of \cite{BTW} in the case where the ambient space
is the Engel group, see also Section~8.1 of the same paper.}
\end{Rem}
%
%
%
%
%
%
%
%
%
\section{Surfaces in the Engel group}\label{sectsurfaces}
%
%
%
%
%
%

In this section we wish to show the following
\begin{The}\label{neglsurf}
Let $\Sigma$ be a 2-dimensional $C^{1,1}$ smooth submanifold of $\E$.
Let $d$ be the degree of $\Sigma$ and let $\Sigma_d$ be the open
subset of points of degree $d$. Then we have
\begin{equation}
\cS^d\big(\Sigma\sm\Sigma_d\big)=0\,.
\end{equation}
\end{The}
In the sequel, we will use the fact that the degree of a point in a submanifold is invariant under left translations.
\begin{Lem}\label{lemcomput}
Let $\Sigma$ be a 2-dimensional $C^1$ smooth submanifold of $\E$ and let $x\in\Sigma$.
Then there there exist local coordinates $u$ in a neighbourhood $U$ of $0$ in $\R^2$ 
such that $x^{-1}\Sigma$ around zero is given by
the local parametrization $\phi:U\lra x^{-1}\Sigma$ with $\phi(0)=0$
and we have
\begin{equation}\label{phi_j}
\phi(u)=
\left\{\begin{array}{ll}
(\phi_1(u),\phi_2(u),u_3,u_4) & \mbox{\rm if $d_\Sigma(x)=5$} \\ 
(u_1,\phi_2(u),\phi_3(u),u_4)\;\mbox{\rm or}\;(\phi_1(u),u_2,\phi_3(u),u_4)
 & \mbox{\rm if $d_\Sigma(x)=4$} \\ 
(u_1,\phi_2(u),u_3,\phi_4(u))\;\mbox{\rm or}\;(\phi_1(u),u_2,u_3,\phi_4(u))
 & \mbox{\rm if $d_\Sigma(x)=3$} \\ 
(u_1,u_2,\phi_3(u),\phi_4(u)) & \mbox{\rm if $d_\Sigma(x)=2$}  \end{array}\right.\,,
\end{equation}
where the functions $\phi_j$'s satisfy
\begin{equation}\label{nablaphi_j}
\left\{\begin{array}{ll}
\nabla\phi_3(0)=\big(0,\der_{u_4}\phi_3(0)\big) & \mbox{\rm if $d_\Sigma(x)=4$} \\ 
\nabla\phi_4(0)=(0,0) & \mbox{\rm if $d_\Sigma(x)=3$} \\ 
\nabla\phi_4(0)=(0,0)\;\mbox{\rm and}\;
\nabla\phi_3(0)=(0,0) & \mbox{\rm if $d_\Sigma(x)=2$}  \end{array}\right..
\end{equation}
Furthermore, for small $r>0$, we have
\begin{equation}\label{vol}
\mu_2(D_{x,r}\cap\Sigma)=r^{d_\Sigma(x)}
\int_{\tilde\delta_{1/r}\big(\phi^{-1}(D_r)\big)}
J\phi(\tilde\delta_ru)\;du\,,
\end{equation}
where $J\phi(x)=\sqrt{\det\big(\lan \der_{u_i}\phi,\der_{u_j}\phi\ran\big)_{i,j=1,\ldots,4}}$
is the Riemannian Jacobian of $\phi$ with respect to the fixed left invariant
Riemannian metric on $\E$. The induced dilations on coordinates $u$ are
defined as $\tilde\delta_r(u)=(r^{d_i}u_i,r^{d_j}u_j)$, where $d_1=d_2=1$,
$d_3=2$ and $d_4=3$.
\end{Lem}
{\sc Proof}.
The proof of \eqref{phi_j} simply follows from the implicit function theorem
and \eqref{formula}.
For example, let us consider the case $d_\Sigma(x)=5$ and let
$\psi$ any local parametrization of $x^{-1}\Sigma$ around the origin.
Then applying \eqref{formula} to $\psi$ at the origin, we must have
$\psi_v^{34}(0)\neq0$, since $\psi(0)=0$. Then the mapping
\[
(v_1,v_2)\lra\big(\psi^3(v),\psi^4(v)\big)
\]
is invertible around the origin and one can take the new coordinates
\[
(u_3,u_4)=(\psi^3(v),\psi^4(v)\big).
\]
The remainig cases proceed in similar way.
Now, if we apply \eqref{formula} to the parametrization $\phi$ at the origin,
having one of the forms given by  \eqref{phi_j}, then a simple computation leads
us to \eqref{nablaphi_j}.
As an example, let us consider the case $d_\Sigma(x)=4$ 
and assume for instance that $\phi(u)=\big(u_1,\phi_2(u),\phi_3(u),u_4\big)$, according to 
the second formula of \eqref{phi_j}.
Then we have
\begin{eqnarray*}
&&\deg\big(\der_{u_1}\phi\wedge\der_{u_4}\phi(0)\big)\\
&&=\deg\left(\big(X_1+\der_{u_1}\phi_2(0)X_2+\der_{u_1}\phi_3(0)X_3\big)\wedge
\big(\der_{u_4}\phi_2(0)X_2+\der_{u_4}\phi_3(0)X_3+X_4\big)\right)\\
&&=\deg\left(X_1\wedge X_4+\der_{u_1}\phi_2(0)X_2\wedge X_4+
\der_{u_1}\phi_3(0)X_3\wedge X_4\right)\,.
\end{eqnarray*}
Since $\deg(X_3\wedge X_4)=5$, then $\der_{u_1}\phi_3(0)=0$ and
the first formula of \eqref{nablaphi_j} follows.
The other cases are achieved in the same way.
The left invariance of the Riemannian surface measure gives
\[
\mu_2\big(\Sigma\cap D_{x,r}\big)=\mu_2\big(x^{-1}\Sigma\cap D_r\big)
=\int_{\phi^{-1}(D_r)}J\phi(u)\;du.
\]
The change of variable $\tilde u=\tilde\delta_r(u)=(r^{d_i}u_i,r^{d_j}u_j)$
and the fact that $d_\Sigma(x)=d_i+d_j$ lead us to formula \eqref{vol}.
\qed

%
%
%
%
%
%

\vskip.2cm
In the sequel, the following box
\[
\Bx_r=[-r,r]^2\times[-r^{2},r^{2}]\times[-r^{3},r^{3}]
\]
will be useful. In fact, by homogeneity there exists $\lambda>0$ such that
\begin{equation}\label{box}
\Bx_{\lambda r}\subset D_r\subset \Bx_{r/\lambda}\quad
\mbox{for every}\quad r>0.
\end{equation}
%
%
%
%
%
%
\begin{Pro}\label{deg45deg2}
Let $\Sigma$ be a 2-dimensional $C^{1,1}$ smooth submanifold of $\E$
and assume that $d(\Sigma)\geq4$ and $d_\Sigma(x)=2$. Then we have
\begin{equation*}
\lim_{r\to0}\dfrac{\mu_2(\Sigma\cap D_{x,r})}{r^{d(\Sigma)}}=+\infty.
\end{equation*}
\end{Pro}
{\sc Proof.} 
We use the coordinates given by Lemma~\ref{lemcomput} and apply 
formulae \eqref{phi_j} and \eqref{nablaphi_j}.
Then we have a constant $c>0$ such that 
$|\phi_3(u)|\leq c|u|^2$ and $|\phi_4(u)|\leq c|u|^2$ for $u$ small.
We have
\begin{equation}
\tilde\delta_{\frac{1}{r}}(\phi^{-1}(\Bx_{\lambda r}))
=\tilde\delta_{1/r}\left\{(u_1,u_2)\!: |u_1|\leq\lambda r,\, |u_2|\leq\lambda r, \,|\phi_3(u)|\leq(\lambda r)^2, \,|\phi_4(u)|\leq(\lambda r)^3\right\},\nonumber
\end{equation}
that can be written as follows
\begin{eqnarray*}
\tilde\delta_{\frac{1}{r}}(\phi^{-1}(\Bx_{\lambda r}))&=&\left\{(x_1,x_2)\!:
\frac{|x_1|}{\lambda}\leq1, \frac{|x_2|}{\lambda}\leq1,
\frac{|\phi_3(rx_1,rx_2)|}{(\lambda r)^2}\leq1,\frac{|\phi_4(rx_1,rx_2)|}
{(\lambda r)^3}\leq1\right\}\\
&\supset&\left\{(x_1,x_2)\!: |x_1|\leq\lambda,\, |x_2|\leq\lambda,
|x|\leq\frac{\lambda}{\sqrt{c}}, |x|\leq\frac{\lambda^{3/2}}{\sqrt{c}}\sqrt{r}\right\}.
\end{eqnarray*}
Taking into account \eqref{vol} and \eqref{box}, it follows that
\[
\dfrac{\mu_2(\Sigma\cap D_{x,r})}{r^{d(\Sigma)}}\geq\frac{\pi\,\lambda^3}{c}
\,\frac{J\phi(0)}{2}\,r^{1-d(\Sigma)+d_\Sigma(x)}=\frac{\pi\,\lambda^3}{c}\,
\frac{J\phi(0)}{2}\,r^{3-d(\Sigma)}\lra+\infty
\quad\mbox{as}\quad r\to 0.
\]
\qed
%
%
%
%
%
%
%
%
\begin{Pro}\label{deg45deg3}
Let $\Sigma$ be a 2-dimensional $C^{1,1}$ smooth submanifold of $\E$
and assume that $d(\Sigma)\geq4$ and $d_\Sigma(x)=3$. Then we have
\begin{equation*}
\lim_{r\to0}\dfrac{\mu_2(\Sigma\cap D_{x,r})}{r^{d(\Sigma)}}=+\infty.
\end{equation*}
\end{Pro}
{\sc Proof}.
From Lemma~\ref{lemcomput}, applying \eqref{phi_j} and \eqref{nablaphi_j},
we get $\phi(u)$ that parametrizes a neighbourhood of 0 in $x^{-1}\Sigma$,
it has the two possible forms
\begin{equation}\label{phideg3}
\phi(u)=(u_1,\phi_2(u),u_3,\phi_4(u))\qquad\mbox{\rm or}\qquad
\phi(u)=(\phi_1(u),u_2,u_3,\phi_4(u))
\end{equation}
and there exists a constant $c>0$ such that 
either $|\phi_2(u)|\leq c|u|$ or $|\phi_1(u)|\leq c|u|$
and also $|\phi_4(u)|\leq c|u|^2$ for $u$ small.
Assume for instance that $\phi(u)=(u_1,\phi_2(u),u_3,\phi_4(u))$.
Of course, the proof is the same in the case $\phi$ assumes
the other form in \eqref{phideg3}. 
The set $\tilde\delta_{\frac{1}{r}}(\phi^{-1}(\Bx_{\lambda r}))$
coincides with 
\begin{eqnarray*}
\left\{(x_1,x_3)\!:\frac{|x_1|}{\lambda}\leq1, \frac{|x_3|}{\lambda^2}\leq1,
\frac{|\phi_2(rx_1,r^2x_3)|}{\lambda r}\leq1,\frac{|\phi_4(rx_1,r^2x_3)|}
{(\lambda r)^3}\leq1\right\}
\end{eqnarray*}
that contains the subset
\begin{eqnarray*}
S_r=\left\{(x_1,x_3)\!: |x_1|\leq\lambda,\, |x_3|\leq\lambda^2,
|(x_1,r x_3)|\leq\frac{\lambda}{c},|(x_1,r x_3)|\leq\sqrt{\frac{\lambda^3 r}{c}}\right\}
\end{eqnarray*}
By the change of variable $x_3'=rx_3$ we get
\begin{equation*}
\cL^2(S_r)=\frac{1}{r}\;
\cL^2\left(\left\{(x_1,x_3')\!: |x_1|\leq\lambda,\, |x_3'|\leq r\lambda^2,
|(x_1,x_3')|\leq\frac{\lambda}{c},|(x_1,x_3')|\leq\sqrt{\frac{\lambda^3 r}{c}}\right\}
\right).
\end{equation*}
For $r>0$ small, it follows that
\begin{eqnarray}\label{srsqest}
\cL^2(S_r)&\geq&\frac{1}{r}\;
\cL^2\left(\left\{(x_1,x_3)\!:  |x_3|\leq r\lambda^2,
\max\{|x_1|,|x_3|\}\leq\sqrt{\frac{\lambda^3 r}{2c}}\right\}\right)\\
&=& \frac{1}{r}\;
\cL^2\left(\left\{(x_1,x_3)\!:  |x_3|\leq r\lambda^2,
|x_1|\leq\sqrt{\frac{\lambda^3 r}{2c}}\right\}\right)\nonumber\\
&=& \frac{2^{3/2}\,\lambda^{7/2}}{\sqrt{c}}\;\sqrt{r}\,.\nonumber
\end{eqnarray}
Finally, taking into account \eqref{srsqest}, \eqref{vol} and \eqref{box},
it follows that
\[
\dfrac{\mu_2(\Sigma\cap D_{x,r})}{r^{d(\Sigma)}}\geq
\frac{\cL^2(S_r)}{r^{d(\Sigma)-3}}\,\frac{J\phi(0)}{2}\lra+\infty
\quad\mbox{as}\quad r\to 0.
\]
\qed
%
%
%
%
%
%
%
%
\begin{Pro}\label{deg3deg2}
Let $\Sigma$ be a 2-dimensional $C^{1,1}$ smooth submanifold of $\E$
and assume that $d(\Sigma)=3$ and $d_\Sigma(x)=2$. Then we have
\begin{equation*}
\lim_{r\to0}\dfrac{\mu_2(\Sigma\cap D_{x,r})}{r^3}=+\infty.
\end{equation*}
\end{Pro}
{\sc Proof.}
From Lemma~\ref{lemcomput}, applying \eqref{phi_j} and \eqref{nablaphi_j},
we get
\begin{equation}\label{phideg2}
\phi(u)=\big(u_1,u_2,\phi_3(u),\phi_4(u)\big)
\end{equation}
that parametrizes a neighbourhood of 0 in $x^{-1}\Sigma$ and
there exists a constant $c>0$ such that 
$|\phi_3(u)|\leq c|u|^2$ and $|\phi_4(u)|\leq c|u|^2$ for $u$ small.
Notice that these estimates do not suffice to obtain our claim.
We have to exploit the assumption that $\Sigma$ has degree three.
In fact, from (\ref{formula}) it follows that
\begin{eqnarray}\label{deg3pdes}
\left\{\begin{array}{l}
\phi_u^{14}-\phi_1\,\phi_u^{13}+\frac{(\phi_1)^2}{2}\,
\phi_u^{12}=0\\
\phi_u^{24}-\phi_1\,\phi_u^{23}=0\\
\phi_u^{34}+\frac{(\phi_1)^2}{2}\phi_u^{23}-\phi_1\phi_u^{24}=0
\end{array}\right..
\end{eqnarray}
Thus, these equations hold for \eqref{phideg2} and in particular
the first two ones yield
\[
\left\{\begin{array}{l}
\partial_{u_2}\phi_4 -u_1\,\partial_{u_2}\phi_3+\frac{(u_1)^2}{2}=0\\
-\partial_{u_1}\phi_4+u_1\partial_{u_1}\phi_3=0\end{array}\right..
\]
We have proved that 
\[
\nabla\phi_4=\left(u_1\partial_{u_1}\phi_3,u_1\,\partial_{u_2}\phi_3
-\frac{(u_1)^2}{2}\right),
\]
hence the $C^{1,1}$ regularity of $\ph$ yieds a constant $c_1>0$ such that $|\phi_4(u)|\leq c_1|u|^3$ for $u$ small. Now, we argue as 
in Proposition~\ref{deg45deg2}, taking into account this better estimate on the
order of vanishing of $\phi_4$. In this case, we have
\begin{eqnarray*}
\tilde\delta_{\frac{1}{r}}(\phi^{-1}(\Bx_{\lambda r}))&=&\left\{(x_1,x_2)\!:
\frac{|x_1|}{\lambda}\leq1, \frac{|x_2|}{\lambda}\leq1,
\frac{|\phi_3(rx_1,rx_2)|}{(\lambda r)^2}\leq1,\frac{|\phi_4(rx_1,rx_2)|}
{(\lambda r)^3}\leq1\right\}\\
&\supset&\left\{(x_1,x_2)\!: |x_1|\leq\lambda,\, |x_2|\leq\lambda,
|x|\leq\frac{\lambda}{\sqrt{c}}, |x|\leq\frac{\lambda}{\sqrt[3]{c_1}}\right\}=S,
\end{eqnarray*}
Then taking into account \eqref{vol} and \eqref{box}, it follows that
\[
\dfrac{\mu_2(\Sigma\cap D_{x,r})}{r^3}\geq
\frac{\cL^2(S)}{r}\,\frac{J\phi(0)}{2}\lra+\infty
\quad\mbox{as}\quad r\to 0.
\]
\qed
\vskip.2truecm
\vskip.2truecm
\noindent
{\sc Proof of Theorem~\ref{neglsurf}}.
Notice that 2-dimensional $C^{1,1}$ submanifolds of degree 2 would
be tangent everywhere to the horizontal subbundle, hence they cannot exist. 
Thus, the possible degrees of $\Sigma$ are 5,4 and 3.
If $d(\Sigma)=5$, then Theorem~2.16 of \cite{Mag5} applies, since
the Hausdorff dimension $Q$ of $\E$ is 7, the codimension of $\Sigma$ is 2
and the set $C(\Sigma)$ in \cite{Mag5} coincides with $\Sigma\sm\Sigma_5$.
Then we have $\cS^5(\Sigma\sm\Sigma_5)=0$. 
If $d(\Sigma)=4$, then Propositions~\ref{deg45deg2} and \ref{deg45deg3}
imply that
\begin{equation}
\lim_{r\to0}\dfrac{\mu_2(\Sigma\cap D_{x,r})}{r^4}=+\infty,
\end{equation}
whenever $x\in\{z\in\Sigma\mid d_\Sigma(z)<4\}=\Sigma\sm\Sigma_4$.
Applying Lemma~\ref{gmta}, we get $\cS^4\big(\Sigma\sm\Sigma_4\big)=0$.
If $d(\Sigma)=3$, then Propositions~\ref{deg3deg2}
implies that
\begin{equation}
\lim_{r\to0}\dfrac{\mu_2(\Sigma\cap D_{x,r})}{r^3}=+\infty,
\end{equation}
whenever $x\in\{z\in\Sigma\mid d_\Sigma(z)<3\}=\Sigma\sm\Sigma_3$.
Then Lemma~\ref{gmta} yields $\cS^3\big(\Sigma\sm\Sigma_3\big)=0$.
This ends the proof.
\qed

%
%
%
%
%
%
%
%
\section{Curves in the Engel group}\label{sectcurves}
%
%
%
%
%
%
%
%

In this section we wish to show the following
\begin{The}\label{neglcurve}
Let $\Sigma$ be a 1-dimensional $C^{1,1}$ smooth submanifold of $\E$.
Let $d$ be the degree of $\Sigma$ and let $\Sigma_d$ be the open
subset of points of degree $d$. Then we have
\begin{equation}
\cS^d\big(\Sigma\sm\Sigma_d\big)=0\,.
\end{equation}
\end{The}
Let $I$ be an open set of $\R$ and let $\phi:I\lra\R^4$ be a $C^1$ immersion.
From computations in Section~4 of \cite{MagVit}, we have
\begin{equation}\label{curvephi}
\dot\phi=\dot\phi_1X_1+\dot\phi_2X_2+\big(\dot\phi_3-\phi_1\dot\phi_2\big)X_3
+\left(\dot\phi_4-\phi_1\dot\phi_3+\frac{(\phi_1)^2}{2}\dot\phi_2\right)X_4\,.
\end{equation}
\begin{Lem}\label{lemcomputcurve}
Let $\Sigma$ be a 1-dimensional $C^1$ smooth submanifold of $\E$ and let $x\in\Sigma$.
Then there exists a neighbourhood $I$ of $0$ in $\R$ 
such that $x^{-1}\Sigma$ around zero is given by
the local parametrization $\phi:I\lra x^{-1}\Sigma$ with $\phi(0)=0$
and we have
\begin{equation}\label{phit}
\phi=
\left\{\begin{array}{ll}
(\phi_1(t_4),\phi_2(t_4),\phi_3(t_4),t_4) & \mbox{\rm if $d_\Sigma(x)=3$} \\ 
(\phi_1(t_3),\phi_2(t_3),t_3,\phi_4(t_3)) & \mbox{\rm if $d_\Sigma(x)=2$} \\ 
(t_1,\phi_2(t_1),\phi_3(t_1),\phi_4(t_1))\;\mbox{\rm or}\;(\phi_1(t_2),t_2,
\phi_3(t_2),\phi_4(t_2)) & \mbox{\rm if $d_\Sigma(x)=1$}
\end{array}\right.\,,
\end{equation}
where the functions $\phi_j$'s satisfy
\begin{equation}\label{dotphit}
\left\{\begin{array}{ll}
\dot\phi_4(0)=0 & \mbox{\rm if $d_\Sigma(x)=2$} \\ 
\dot\phi_4(0)=\dot\phi_3(0)=0 & \mbox{\rm if $d_\Sigma(x)=1$} 
\end{array}\right..
\end{equation}
Furthermore, for small $r>0$, we have
\begin{equation}\label{length}
\mu_1(D_{x,r}\cap\Sigma)=r^{d_\Sigma(x)}
\int_{\tilde\delta_{1/r}\big(\phi^{-1}(D_r)\big)}
J\phi(\tilde\delta_rt)\;dt\,,
\end{equation}
where $J\phi(t)=\sqrt{\lan\dot\phi(t),\dot\phi(t)\ran}$
and $\lan,\ran$ denotes the fixed left invariant Riemannian metric.
The induced dilation is $\tilde\delta_r(t_j)=r^{d_j}t_j$,
where $d_1=d_2=1$, $d_3=2$ and $d_4=3$.
\end{Lem}
{\sc Proof}.
Taking into account the notion of pointwise degree, the proof of \eqref{phit}
simply follows from the implicit function theorem and formula \eqref{curvephi}.
As in the proof of Lemma~\ref{lemcomput}, we consider just one case as an example.
The other ones, can be simply obtained repeating this argument.
For instance, if $d_\Sigma(x)=3$, then $\dot\psi_4(0)\neq0$,
where $\psi$ is a parametrization of $x^{-1}\Sigma$ around the origin.
Setting $t_4=\psi(s)$, we have proved \eqref{phit} in the case $d_\Sigma(x)=3$.
Now, if we apply \eqref{curvephi} to the parametrization $\phi$,
at the origin, then a simple computation gives \eqref{dotphit}.
In fact, one has applied \eqref{curvephi} to all forms of $\phi$
listed in \eqref{phit}.
The left invariance of the Riemannian surface measure gives
\[
\mu_1\big(\Sigma\cap D_{x,r}\big)=\mu_1\big(x^{-1}\Sigma\cap D_r\big)
=\int_{\phi^{-1}(D_r)}J\phi(t)\;dt.
\]
The change of variable $\tilde t_j=r^{d_j}t_j$
and the fact that $d_\Sigma(x)=d_j$ lead us to formula \eqref{length}.
\qed
%
%
%
%
%
%
\begin{Pro}\label{deg23deg1}
Let $\Sigma$ be a 1-dimensional $C^{1,1}$ smooth submanifold of $\E$
and assume that $d(\Sigma)\geq2$ and $d_\Sigma(x)=1$. Then we have
\begin{equation*}
\lim_{r\to0}\dfrac{\mu_1(\Sigma\cap D_{x,r})}{r^{d(\Sigma)}}=+\infty.
\end{equation*}
\end{Pro}
{\sc Proof.} 
By Lemma~\ref{lemcomputcurve}, applying formulae \eqref{phit} and \eqref{dotphit},
we get a constant $c>0$ such that 
$\max\{|\phi_3(t)|,|\phi_4(t)|\}\leq c|t|^2$ 
and $\phi_2(t)\leq c|t|$ for $t$ small.
Moreover, $x^{-1}\Sigma$ can be locally parametrized either by
$(t_1,\phi_2(t_1),\phi_3(t_1),\phi_4(t_1))$ or
$\phi_1(t_2),t_2,\phi_3(t_2),\phi_4(t_2))$.
Clearly, both cases have the same proof, hence assume for instance that
\[
\phi(t_1)=(t_1,\phi_2(t_1),\phi_3(t_1),\phi_4(t_1))
\]
locally parametrizes $x^{-1}\Sigma$.
Taking into account \eqref{box}, we consider
\begin{equation}
\tilde\delta_{\frac{1}{r}}(\phi^{-1}(\Bx_{\lambda r}))
=\tilde\delta_{1/r}\left\{t_1\!:
\frac{|t_1|}{\lambda r}\leq1,\, \frac{|\phi_2(t_1)|}{\lambda r}\leq1,
\,\frac{|\phi_3(t_1)|}{(\lambda r)^2}\leq1,
\,\frac{|\phi_4(t_1)|}{(\lambda r)^3}\leq1 \right\}.\nonumber
\end{equation}
This set can be written as
\begin{eqnarray*}
\left\{\tau :
\frac{|\tau|}{\lambda}\leq1,\, \frac{|\phi_2(r\tau)|}{\lambda r}\leq1,
\,\frac{|\phi_3(r\tau)|}{(\lambda r)^2}\leq1,
\,\frac{|\phi_4(r\tau)|}{(\lambda r)^3}\leq1 \right\}.\nonumber
\end{eqnarray*}
For $r>0$ small, the previous set contains
\begin{eqnarray*}
S_r=\left\{\tau : |\tau|\leq\frac{\lambda^{3/2}}{\sqrt{c}}\sqrt{r}\right\}.\nonumber
\end{eqnarray*}

Taking into account \eqref{length} and \eqref{box}, it follows that
\[
\dfrac{\mu_1(\Sigma\cap D_{x,r})}{r^{d(\Sigma)}}\geq
\frac{\cL^1(S_r)}{r^{d(\Sigma)-1}}
\,\frac{J\phi(0)}{2}\lra+\infty
\quad\mbox{as}\quad r\to 0.
\]
\qed
%
%
%
%
%
%
\begin{Pro}\label{deg3deg2curve}
Let $\Sigma$ be a 1-dimensional $C^{1,1}$ smooth submanifold of $\E$
and assume that $d(\Sigma)=3$ and $d_\Sigma(x)=2$. Then we have
\begin{equation*}
\lim_{r\to0}\dfrac{\mu_1(\Sigma\cap D_{x,r})}{r^3}=+\infty.
\end{equation*}
\end{Pro}
{\sc Proof.} 
By Lemma~\ref{lemcomputcurve}, applying formulae \eqref{phit} and \eqref{dotphit},
we get a constant $c>0$ such that 
$|\phi_4(t)|\leq c|t|^2$ and $\max\{|\phi_1(t)|,|\phi_2(t)|\}\leq c|t|$ for $t$ small.
Moreover, $x^{-1}\Sigma$ can be locally parametrized by
\[
(\phi_1(t_3),\phi_2(t_3),t_3,\phi_4(t_3)).
\]
Taking into account the previous estimates on $\phi_j$,
for $r>0$ small we have
\begin{eqnarray*}
\tilde\delta_{\frac{1}{r}}(\phi^{-1}(\Bx_{\lambda r}))
&=&\tilde\delta_{1/r}\left\{t_1\!:
\frac{|\phi_1(t_3)|}{\lambda r}\leq1,\, \frac{|\phi_2(t_3)|}{\lambda r}\leq1,
\,\frac{|t_3|}{(\lambda r)^2}\leq1,
\,\frac{|\phi_4(t_3)|}{(\lambda r)^3}\leq1 \right\}\nonumber\\
&=&\left\{\tau : |\tau|\leq\lambda\right\}=S.\nonumber
\end{eqnarray*}
Thus, from \eqref{length} and \eqref{box}, it follows that
\[
\dfrac{\mu_1(\Sigma\cap D_{x,r})}{r^3}\geq
\frac{\cL^1(S)}{r}
\,\frac{J\phi(0)}{2}\lra+\infty
\quad\mbox{as}\quad r\to 0.
\]
\qed
\vskip.2truecm
\noindent
{\sc Proof of Theorem~\ref{neglcurve}}.
If $d(\Sigma)=3$, then Propositions~\ref{deg23deg1} and \ref{deg3deg2}
imply that
\begin{equation}
\lim_{r\to0}\dfrac{\mu_1(\Sigma\cap D_{x,r})}{r^3}=+\infty,
\end{equation}
whenever $x\in\{z\in\Sigma\mid d_\Sigma(z)<4\}=\Sigma\sm\Sigma_4$.
Thus, Lemma~\ref{gmta} yields $\cS^3\big(\Sigma\sm\Sigma_3\big)=0$.
If $d(\Sigma)=2$, then Proposition~\ref{deg23deg1}
implies that
\begin{equation}
\lim_{r\to0}\dfrac{\mu_2(\Sigma\cap D_{x,r})}{r^2}=+\infty,
\end{equation}
whenever $x\in\{z\in\Sigma\mid d_\Sigma(z)<2\}=\Sigma\sm\Sigma_2$
and this gives $\cS^2\big(\Sigma\sm\Sigma_2\big)=0$.
\qed

\end{document}